\documentclass[11pt]{article}
\usepackage{amsmath,amssymb,theorem,amstext,amsgen,amsbsy,amsopn,amsfonts,graphicx,cases}
\usepackage{graphicx}
\usepackage{subfigure}
\usepackage{psfrag}
\usepackage{enumerate}
\usepackage{color}
\usepackage{epstopdf}

\newtheorem{thm}{Theorem}[section]

\newtheorem{obs}[thm]{Observation}
\newtheorem{cor}[thm]{Corollary}
\newtheorem{prop}[thm]{Proposition}
\theorembodyfont{\rmfamily}

\def\pf{\bigskip\noindent {\bf Proof.}}

\def\dfn#1{{\sl #1}}

\def\less{\backslash}

\def\qed{ \hfill\vrule height3pt width6pt depth2pt}

\def\pf{\bigskip\noindent {\bf Proof.}~~}

\title{Extremal   Theta-free planar graphs}

\author{Yongxin Lan and  Yongtang Shi  \\ 
 Center for Combinatorics and LPMC\\
Nankai University, Tianjin 300071, China\\
and\\
  Zi-Xia Song\thanks{Corresponding author.    E-mail address: Zixia.Song@ucf.edu}\\
 Department  of Mathematics\\
 University of Central Florida, Orlando, FL 32816, USA
}

\begin{document}
\maketitle
\begin{abstract}
 Given a family $\mathcal{F}$, a graph is   \dfn{$\mathcal{F}$-free}  if it does not contain any graph in $\mathcal{F}$ as a subgraph. We   study the topic of ``extremal" planar graphs initiated by Dowden [J. Graph Theory  83 (2016) 213--230], that is,   how many edges
can an $\mathcal{F}$-free planar graph on $n$ vertices have?   We define $ex_{_\mathcal{P}}(n,\mathcal{F})$ to be  the maximum number of edges in an $\mathcal{F}$-free planar graph on $n $ vertices.
  Dowden   obtained the tight bounds  $ex_{_\mathcal{P}}(n,C_4)\leq15(n-2)/7$ for all $n\geq4$ and $ex_{_\mathcal{P}}(n,C_5)\leq(12n-33)/5$ for all $n\geq11$.
 In this paper,  we continue to promote the idea of determining  $ex_{_\mathcal{P}}(n,\mathcal{F})$ for certain classes $\mathcal{F}$.  Let $\Theta_k$ denote  the family of  Theta graphs on $k\ge4$ vertices, that is,  graphs obtained from a cycle $C_k$ by adding  an additional edge joining two non-consecutive vertices. The study of $ex_{_\mathcal{P}}(n,\Theta_4)$ was suggested by Dowden.
 We show that $ex_{_\mathcal{P}}(n,\Theta_4)\leq12(n-2)/5$ for all $n\geq 4$,  $ex_{_\mathcal{P}}(n,\Theta_5)\leq5(n-2)/2$ for all $n\ge5$, and then demonstrate that these bounds are tight, in the sense that there are infinitely many values of $n$ for which they are attained exactly.
 We also prove that  $ex_{_\mathcal{P}}(n,C_6)\le ex_{_\mathcal{P}}(n,\Theta_6)\le 18(n-2)/7$ for all $n\ge6$.
\end{abstract}


\baselineskip 17pt
\section{Introduction}

All graphs considered in this paper are finite and simple.  We use  $C_k$ to denote the cycle on $k$ vertices. Let $\mathcal{F}$ be a family of graphs. A graph is $\mathcal{F}$-free if it does not contain any graph in $\mathcal{F}$ as a subgraph. When $\mathcal{F}=\{F\}$ we write $F$-free.  One of the best known results in extremal graph theory is Tur\'an's Theorem~\cite{1941Turan}, which
gives the maximum number of edges that a $K_k$-free graph on $n$ vertices can have. The celebrated  Erd\H{o}s-Stone Theorem~\cite{1946ErdosStone}  then extends this to the case when $K_k$ is
replaced by an arbitrary graph $H$, showing that the maximum number of edges possible
is $(1+o(1)) (\frac{\chi(H)-2}{\chi(H)-1})  {n\choose 2} $, where $\chi(H)$ denotes the chromatic number of $H$. This latter result has been called the ``fundamental theorem of extremal graph theory''~\cite{2013Bollobas}.  Tur\'an-type problems when host graphs are hypergraphs are notoriously difficult.  A large quantity of work in this area has been carried out  in determining the maximum number of edges in a $k$-uniform  hypergraph  on $n$ vertices without containing  $k$-uniform linear  paths and  cycles (see, for example, \cite{FurediJiang,FurediJiangSeiver,Kostochka1}). Surveys on Tur\'an-type problems of graphs and hypergraphs can be found in \cite{Furedi} and \cite{Keevash}.\medskip

Recently, Dowden~\cite{Dowden} initiated the study of Tur\'an-type problems when host graphs are planar graphs, i.e.,  how many edges can an $\mathcal{F}$-free  planar graph on $n$ vertices have?  The \dfn{planar Tur\'an number of  $\mathcal{F}$}, denoted  $ex_{_\mathcal{P}}(n,\mathcal{F})$, is  the maximum number of edges in an $\mathcal{F}$-free  planar graph on $n$ vertices.   When $\mathcal{F}=\{F\}$ we write $ex_{_\mathcal{P}}(n,F)$. Dowden~\cite{Dowden} observed that it is straightforward to determine the exact values of $ex_{_\mathcal{P}}(n,H)$ when $H$ is a complete graph or non-planar graph; he also obtained the tight bounds  $ex_{_\mathcal{P}}(n,C_4)\leq15(n-2)/7$ for all $n\geq4$ and $ex_{_\mathcal{P}}(n,C_5)\leq(12n-33)/5$ for all $n\geq11$.
 Recently, Lan, Shi and Song observed in \cite{2017LSZ} that  planar Tur\'an numbers are closely related to  planar anti-Ramsey numbers. The \dfn{planar anti-Ramsey number of $\mathcal{F}$}, denoted
 $ar_{_\mathcal{P}}(n,\mathcal{F})$,     is the  maximum number  of colors in an edge-coloring of  a plane triangulation $T $  (which is not $\mathcal{F}$-free) on $n$ vertices  such that $T$ contains   no    rainbow copy of any $F\in\mathcal{F}$.    When $\mathcal{F}=\{F\}$ we write $ar_{_\mathcal{P}}(n,F)$. The study of planar anti-Ramsey numbers was initiated by Hor\v{n}\'ak,  Jendrol$'$,  Schiermeyer and  Sot\'ak~\cite{HJSS} (under the name of rainbow numbers).  The following   result is observed in \cite{2017LSZ}.

 \begin{prop}[\cite{2017LSZ}]\label{LU} Given a planar graph $H$ and a positive integer $n\ge |H|$,
 $$1+ex_{_\mathcal{P}}(n,\mathcal{H})\le ar_{_\mathcal{P}}(n,H) \le ex_{_\mathcal{P}}(n, H),$$
 where $\mathcal{H}=\{H-e: \, e\in E(H)\}$.
 \end{prop}

In this paper,  we continue to promote the idea of determining  $ex_{_\mathcal{P}}(n,\mathcal{F})$ for certain classes $\mathcal{F}$.  This paper focuses   on    the family  of Theta graphs, where  a graph on at least $4$ vertices is a \dfn{Theta graph} if it can be obtained from a cycle  by adding  an additional edge joining two non-consecutive vertices.  For integer $k\ge4$, let $\Theta_k$ be the family of non-isomorphic Theta graphs on $k$ vertices.   Note  that the only graph in $\Theta_4$ is isomorphic to $K_4$ minus one edge, and $\Theta_5$ has only one graph. By abusing notation, we also use  $\Theta_4$ and $\Theta_5$ to denote the only graph in $\Theta_4$ and $\Theta_5$, respectively.   Note that the study of $ex_{_\mathcal{P}}(n,\Theta_4)$ was suggested by Dowden~\cite{Dowden}.  We  state and prove our main results in Section~\ref{ThetaGraphs}. More recent results on this topic can be found in \cite{2018LSZ}. However,   determining  $ex_{_\mathcal{P}}(n,H)$,   when $H$ is a  planar subcubic graph,  remains wide open. \medskip

We need to introduce more notation.  For a graph $G$, we will use $V(G)$ to denote the vertex set, $E(G)$ the edge set, $|G|$ the number of vertices, $e(G)$ the number of edges and $\delta(G)$ the minimum degree.  For a vertex $x \in V(G)$, we will use $N_G(x)$ to denote the set of vertices in $G$ which are adjacent to $x$.
We define  $d_G(x) = |N_G(x)|$. Given vertex sets $A, B \subseteq V(G)$,
the subgraph of $G$ induced on  $A$, denoted $G[A]$, is the graph with vertex set $A$ and edge set $\{xy \in E(G) : x, y \in A\}$. We denote by $B \less A$ the set $B - A$  and $G \less A$ the subgraph of $G$ induced  on
$V(G) \less A$, respectively.      
Given two isomorphic graphs $G$ and $H$,   we may (with a slight but common abuse of notation) write $G = H$. 
  For any positive integer $k$, let  $[k]:=\{1,2, \ldots, k\}$. \medskip

An $\mathcal{F}$-free  planar graph $G$ on $n$ vertices with the largest possible number of edges is called \dfn{extremal for $n$ and $\mathcal{F}$}. If $\mathcal{F}=\{F\}$, then we simply say $G$ is extremal for $n$ and $F$. Given a plane graph $G$ and integers $i,j\ge3$, an \dfn{$i$-face} in $G$  is  a face of size $i$; and  let:
$E_{i,j}$ denote the set of edges  in $G$ that each  belong  to  one  $i$-face and one  $j$-face (and  belong  to  two  $i$-faces when $i=j$);  $E_i$ denote the set of edges in $G$ that each  belong  to at least one $i$-face; and  $f_i$  denote the number of $i$-faces  in $G$.  Let  $e_{i,j}: =|E_{i,j}|$, $e_i:=|E_i|$,  and   $f:=\sum_if_i$.  Given three positive integers $a, b$ and $c$, we use  $a\equiv b (\rm{mod}\, c)$ to denote $a$ and $b$ have  the  same remainder when   divided by $c$.  We will make use of the following observation.

\begin{obs}\label{fact} Let $G$ be a plane graph on $n\ge3$ vertices  with $e(G)\ge2$.
 For all $i\ge3$,
\begin{enumerate}[(a)]
 \item $e_{i,i}\le e_i\le e(G)$,\vspace{-2mm}
 \item  $if_i=e_i+e_{i,i}$,   \vspace{-2mm}
 \item $\sum_{ i\ge3}{e_i}-\sum_{ 3\le i<j}  e_{i, j}= e(G)$, and \vspace{-1mm}
 \item every face in  $G$  is bounded by a cycle if $G$ is  $2$-connected.
 \end{enumerate}
\end{obs}


\section{Planar Tur\'an numbers of Theta graphs}\label{ThetaGraphs}

   In this section, using the method developed in \cite{Dowden}, we study planar Tur\'an numbers of $\Theta_k$ when  $k\in\{4,5,6\}$. We begin with $\mathcal {F}=\Theta_4$ and prove that  $ex_{_\mathcal{P}}(n,\Theta_4)\leq {12(n-2)}/5$  for all $n\ge4$ and   then demonstrate that this bound is tight, in the sense that there are infinitely many values of $n$ for which it is attained exactly.

\begin{thm}\label{Theta4}
   $ex_{_\mathcal{P}}(n,\Theta_4)\leq  {12(n-2)}/5$ for all $n\ge4$, with  equality when $n\equiv 12 (\rm{mod}\, 20)$.
\end{thm}
\pf  Let $G$ be a $\Theta_4$-free plane graph on $n\geq4$ vertices.  We shall proceed the proof by induction on $n$. 
The statement is trivially true when $n=4$ because any $\Theta_4$-free plane graph on four vertices has at most four edges. So we may assume that  $n\ge5$. Next assume that  there exists a vertex $u\in V(G)$ with  $d_G(u)\leq2$.
By the induction hypothesis, $ e(G\less u)\leq 12(n-3)/5$ and so $e(G)= e(G\less u)+d_G(u)\leq 12(n-3)/5+2<12(n-2)/5$, as desired. So we may assume that $\delta(G)\geq3$. Then each component of $G$ has at least five vertices because $G$ is $\Theta_4$-free.
By the induction hypothesis,  we may further assume that $G$ is connected.  Then $G$ has no face of size at most two because $\delta(G)\geq3$. Hence 

\begin{equation}\tag{1}
2e(G)=\sum_{i\geq3}if_i
\geq 3f_3+4\sum_{i\geq4}f_i
=3f_3+4(f-f_3)
=4f-f_3,
\end{equation}
which implies that $f\leq (2e(G)+f_3)/4$. Note that  $E_{3,3}=\emptyset$ else $G$ would contain $\Theta_4$ as a subgraph, a contradiction. Thus $e_3=3f_3$ by Observation~\ref{fact}(b). This, together with $e_3\le e(G)$ and $f\leq (2e(G)+f_3)/4$, implies that $f\leq 7e(G)/12$.
By Euler's formula, $n-2=e(G)-f\geq 5e(G)/12$. Hence $e(G)\leq 12(n-2)/5$, as desired.\\

\begin{figure}[htbp]
\centering
\subfigure[][\label{G0}]{
\hfill\includegraphics[scale=0.35]{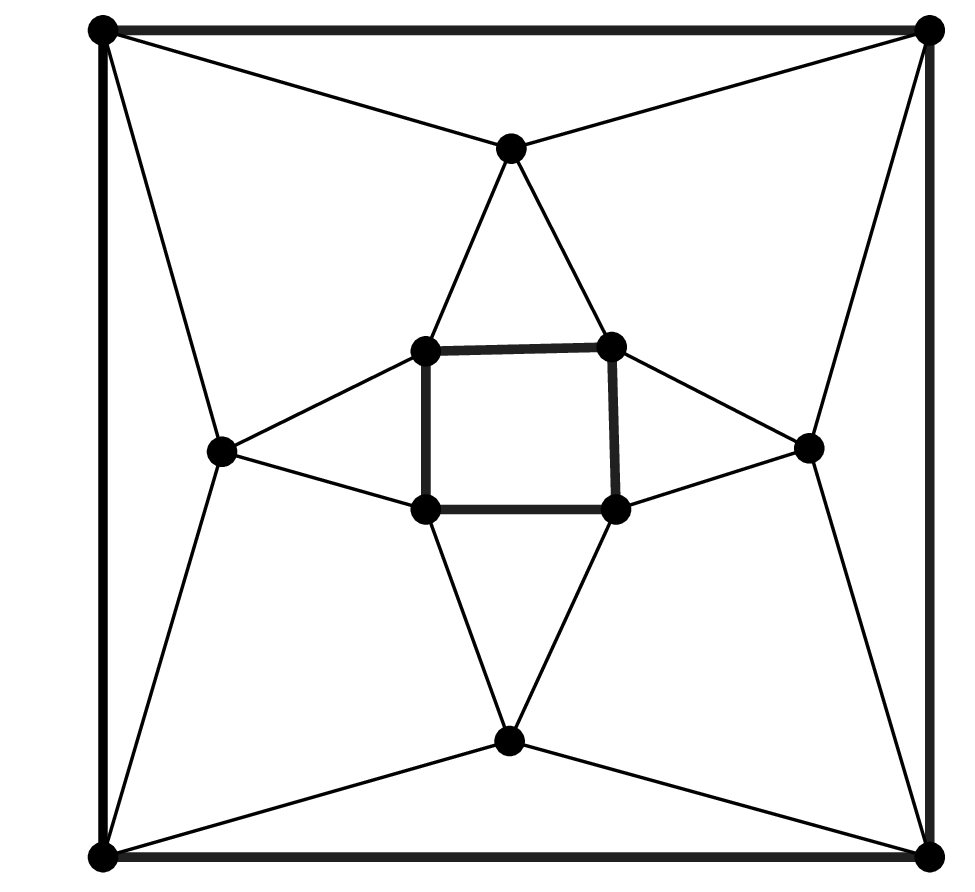}
}
\hskip 1cm
\subfigure[][\label{Gk}]{
\hfill\includegraphics[scale=0.35]{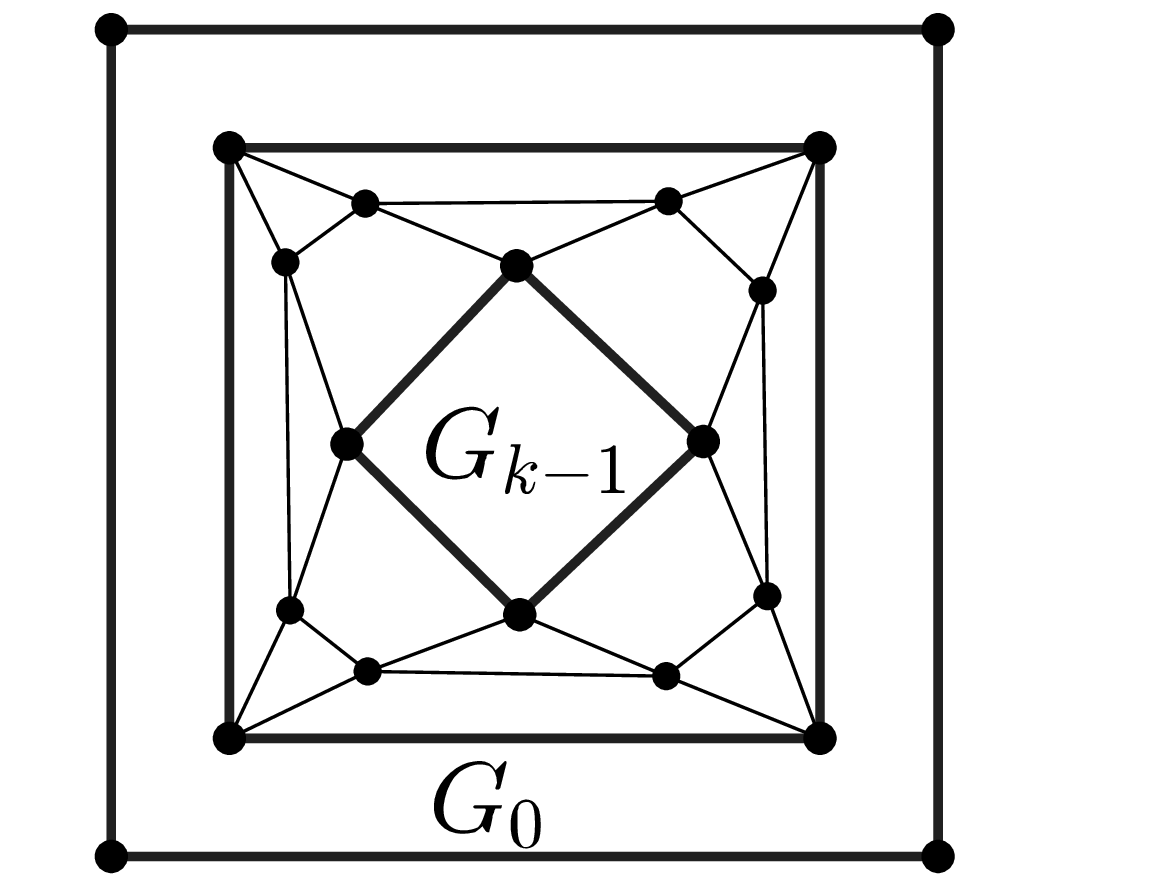}
}
\caption{Construction of $G_k$.}
\end{figure}

From the proof above, we see that   equality in $e(G)\le  12(n-2)/5$ is achieved for $n$  if and only if  equalities  hold both in (1) and in $e_3\le e(G)$. This implies that $e(G)= 12(n-2)/5$ for $n$  if and only if $G$ is a  connected $\Theta_4$-free plane graph on $n$ vertices such that each edge in $G$ belongs to  one $3$-face and one $4$-face. We next construct such an extremal  graph for $n$ and $\Theta_4$. Let $n=20k+12$ for some  integer $k\ge0$.  Let $G_0$ be the graph depicted in Figure~\ref{G0},  we then construct $G_k$ on $n$ vertices recursively  for  all $k\ge1$ via the illustration given in Figure~\ref{Gk}: the entire graph $G_{k-1}$ is placed  into the center quadrangle of Figure~\ref{Gk}, and the entire  $G_0$ is then placed  between  the two given bold quadrangles of Figure~\ref{Gk}  (in such a way that these are identified with the bold quadrangles of  Figure~\ref{G0}).  One can check that $G_k$ is $\Theta_4$-free with  $n=20k+12$ vertices and $ 12(n-2)/5$ edges for all $k\ge0$.  \qed\\

We next prove that  $ex_{_\mathcal{P}}(n,\Theta_5)\leq{5(n-2)}/2$ and   then demonstrate that this bound is tight, in the sense that there are infinitely many values of $n$ for which it is attained exactly.

\begin{thm}\label{Theta5}
   $ex_{_\mathcal{P}}(n,\Theta_5)\leq {5(n-2)}/2$ for all $n\ge5$, with  equality when $n\equiv 50 (\rm{mod}\, 120)$.
\end{thm}

\pf Let $G$ be a $\Theta_5$-free plane graph on $n\geq 5$ vertices. We show by induction on $n$ that $e(G)\le {5(n-2)}/2$. The statement is trivially true when $n=5$ because any $\Theta_5$-free plane graph on five vertices has at most seven edges. So we may assume that $n\ge6$. Next assume that  there exists a vertex $u\in V(G)$ with  $d_G(u)\leq2$. By the induction hypothesis, $ e(G\less u)\leq 5(n-3)/2$ and so  $e(G)= e(G\less u)+d_G(u)\leq 5(n-3)/2+2<5(n-2)/2$, as desired. So we may assume that $\delta(G)\geq3$. Assume next that  $G$ is disconnected. Then each  component of $G$ has either exactly four vertices or at least six vertices because $G$ is $\Theta_5$-free and $\delta(G)\ge 3$. Let $G_1, \dots, G_s, G_{s+1}, \dots, G_{s+t}$ be all components of $G$ such that $|G_1|=\cdots=|G_s|=4$ and $6\le |G_{s+1}|\le\cdots\le |G_{s+t}|$, where $s\ge 0$ and $t\ge0$ are integers with $s+t\ge2$ and $4s+|G_{s+1}|+\cdots+|G_{s+t}|=n$.  Then  $e(G_i)=6$ for all $i\in[s]$ because $\delta(G)\ge 3$,  and  $e(G_j)\le 5(|G_j|-2)/2$ for all $j\in\{s+1, \dots, s+t\}$ by the induction hypothesis. Therefore,
\begin{align*}
e(G)&\le 6s+\frac{5(|G_{s+1}|+\cdots+|G_{s+t}|-2t)}2\\
&=\frac{5(n-2)}2-\frac{(8(s+t)+2t-10)}2\\
&<\frac{5(n-2)}2,
\end{align*}
as desired.  So we may further assume that $G$ is connected.\medskip

Since $G$ is  a connected plane graph  on $n\ge 6$ vertices,  we see that $G$ has no face of size at most two. Hence
\begin{equation}\tag{2}
2e(G)=3f_3+4f_4+\sum_{i\geq 5}if_i
\geq 3f_3+4f_4+5(f-f_3-f_4)
=5f-2f_3-f_4,
\end{equation}
which implies that $f\leq (2e(G)+2f_3+f_4)/5$.
Note that  no $3$-face in $G$ has its  three edges in $E_{3,3}$ because   $G$ is $\Theta_5$-free and $n\ge6$.  It follows that $e_{3,3}\le f_3$. By Observation~\ref{fact}(b),
\begin{equation}\tag{3}
3f_3=e_3+e_{3,3}\leq e_3+f_3 \text{ and so } f_3\leq e_3/2.
\end{equation}

It is worth noting that
 a $4$-face and a $3$-face in G cannot have exactly one edge in common,
 else $G$ would contain $\Theta_5$ as a subgraph. Hence, $E_{3,4}=\emptyset$.  Since $\delta(G)\ge3$, we see that a $4$-face and a $3$-face in $G$    cannot have exactly two  edges in common. Hence  every  $4$-face and  every $3$-face in $G$ have no edge in common,   and so $e_3+e_4\leq e(G)$. By Observation~\ref{fact}(a,b), $e_{4,4}\le e_4$ and  $4f_4=e_4+e_{4,4}$. It follows that
\begin{equation}\tag{4}
4f_4\le 2e_4\le 2(e(G)-e_3) \text{ and so } f_4\leq (e(G)-e_3)/2.
\end{equation}
Now with the last inequalities in (3) and (4), and the fact that  $f\leq (2e(G)+2f_3+f_4)/5$ and $e_3\leq e(G)$, we obtain  $f\leq 3e(G)/5$. By Euler's formula, $n-2=e(G)-f\geq 2e(G)/5$. Hence $e(G)\leq 5(n-2)/2$, as desired.\medskip

\begin{figure}[htbp]
\centering
\subfigure[][\label{G5}]{
\hfill\includegraphics[scale=0.27]{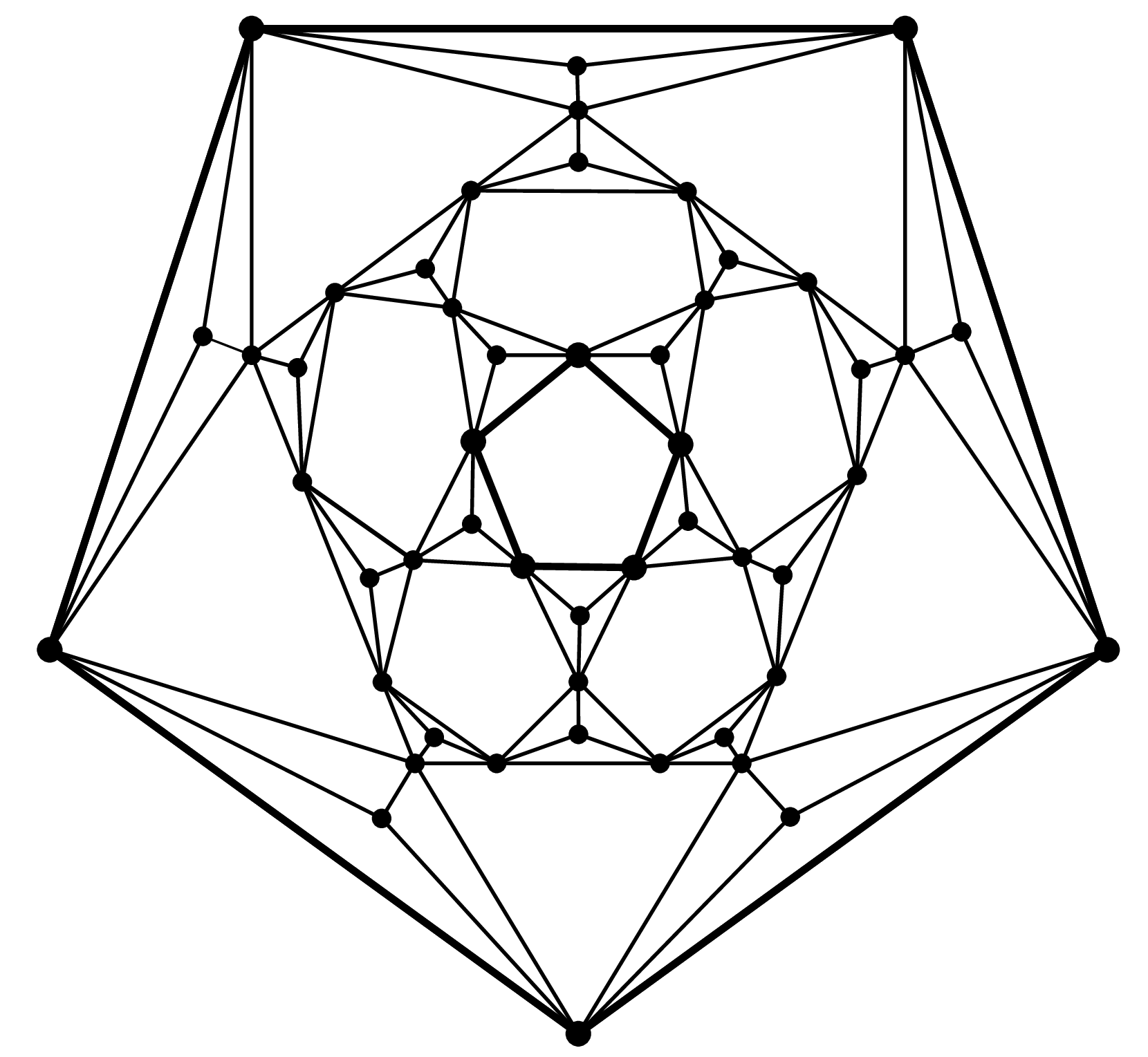}
}
\hskip 0.6cm
\subfigure[][\label{G5k}]{
\hfill\includegraphics[scale=0.27]{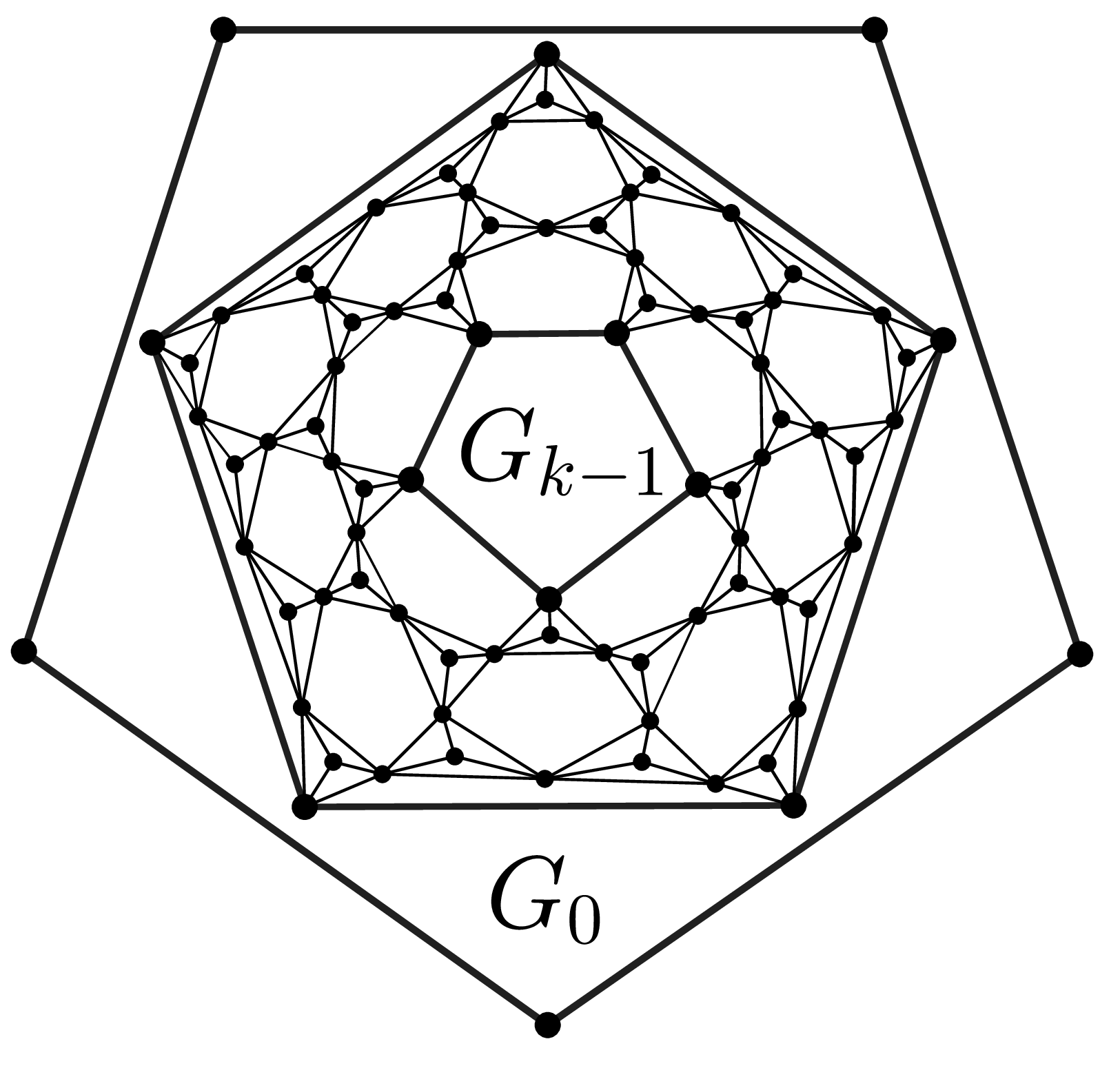}
}
\caption{Construction of $G_k$.}
\end{figure}

From the proof above, we see that  equality in $e(G)\le 5(n-2)/2$ is achieved for $n$  if and only if equalities  hold  in (2), (3) and (4) and in $e_3\le e(G)$. This implies that $e(G)= 5(n-2)/2$  for $n$  if and only if $G$ is a connected $\Theta_5$-free plane graph on $n$ vertices satisfying: each $3$-face in $G$  has exactly two edges in $E_{3,3}$; each edge in $G$ belongs to  either one $3$-face and one $5$-face or two $3$-faces. We next construct such an extremal plane graph for $n$ and $\Theta_5$.  Let $n=120k+50$ for some  integer $k\ge0$. Let $G_0$ be the graph depicted in Figure~\ref{G5}, we then construct $G_k$ of order $n$ recursively  for  all $k\ge1$ via the illustration given in Figure~\ref{G5k}: the entire graph $G_{k-1}$ is placed  into the center pentagon of Figure~\ref{G5k}, and the entire $G_0$ is then placed  between  the  two given bold pentagons of Figure~\ref{G5k}  (in such a way that these are identified with the bold pentagons of  Figure~\ref{G5}).  One can check that $G_k$ is $\Theta_5$-free  with $n=120k+50$ vertices and $ 5(n-2)/2$ edges for all $k\ge0$. \qed\\

\begin{figure}[htbp]
\centering
\includegraphics [scale=0.7] {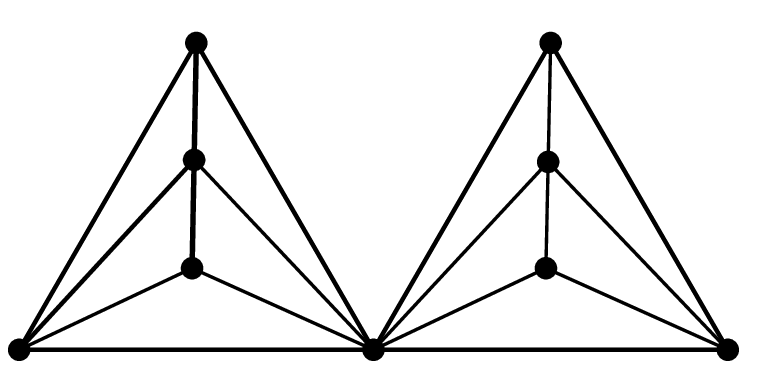}
\caption{A graph   achieving equality in Theorem~\ref{Theta6} and Corollaries~\ref{strict}, \ref{C6} when  $n=9$.}\label{fig3}
\end{figure}

Finally, we prove an upper bound for $ex_{_\mathcal{P}}(n,\Theta_6)$ in Theorem~\ref{Theta6}.  Figure~\ref{fig3}  gives an example of a graph for which equality in Theorem~\ref{Theta6}  is attained when $n=9$. However, we shall  see in Corollary~\ref{strict} that  equality is not possible for larger values of $n$.

\begin{thm}\label{Theta6}
   $ex_{_\mathcal{P}}(n,\Theta_6)\leq 18(n-2)/7$ for all $n\geq6$, with  equality  when $n=9$.
\end{thm}
\begin{pf}
Let $G$ be an  extremal  plane graph  for $\Theta_6$ and $n\geq6$. We shall prove   that $e(G)\le {18(n-2)}/7$    by induction on $n$. The statement is trivially true for $n=6$ because any $\Theta_6$-free plane graph on six vertices has at most ten edges. So we may assume that $n\ge7$.   Next assume that  there exists a vertex $u\in V(G)$ with  $d_G(u)\leq2$. By the induction hypothesis,  $ e(G\less u) \leq 18(n-3)/7$ and so $e(G)= e(G\less u)+d_G(u)\leq 18(n-3)/7+2<18(n-2)/7$, as desired. So we may assume that $\delta(G)\geq3$. Assume next that  $G$ is disconnected.  Then each component of $G$ has exactly four, five or at least six vertices  because $\delta(G)\ge 3$.
 Let $G_1, \dots, G_r, G_{r+1}, \dots, G_{r+s}, G_{r+s+1}, \dots, G_{r+s+t}$ be all components of $G$ such that
\[
|G_1|= \cdots=|G_r|=4, |G_{r+1}|=\cdots= |G_{r+s}| =5,  \text{ and } 6\le |G_{r+s+1}|\le\cdots\le |G_{r+s+t}|,
\]
where $r,s,t\ge 0$ are integers with $r+s+t\ge2$ and $4r+5s+|G_{r+s+1}|+\cdots+|G_{r+s+t}|=n$.  Since $G$ is  an  extremal  plane graph for $\Theta_6$, we see that  $e(G_i)=6$ for all $i\in[r]$ and  $e(G_j)=9$ for all $j \in \{r+1, \ldots, r+s\}$.  By the induction hypothesis,  $e(G_k)\le 18(|G_k|-2)/7$ for all $k\in\{r+s+1, \ldots, r+s+t\}$.   Therefore,
\begin{align*}
e(G)&\le 6r+9s+\frac{18(|G_{r+s+1}|+\cdots+|G_{r+s+t}|-2t)}7\\
&=\frac{18(n-2)}7-\frac{(27(r+s+t)+3r+ 9t -36)}7\\
&<\frac{18(n-2)}7,
\end{align*}
as desired. So we may  assume that $G$ is connected. \medskip

Next assume  that $G$ contains a cut-vertex, say $u$.  Let $H$ be a  smallest component of $G\less u$, and let  $G_1: =G[V(H)\cup \{u\}]$ and $G_2:=G\less V(H)$. Then  $|G_1|\le |G_2|$ and $|G_1|+|G_2|=n+1$.  Since $\delta(G)\geq3$, we see that  $4\leq |G_1|\leq |G_2|$.  Assume first that  $|G_2|\leq 5$. Then $e(G_i)\leq 3|G_i|-6$ for all $i\in\{1,2\}$. Hence,
$e(G)=e(G_1)+e(G_2)\leq 3(|G_1|+|G_2|)-12=3n-9\le18(n-2)/7$ because $n\le9$, with equality when both $G_1$ and $G_2$ are isomorphic to     $K_5$ minus  one edge, and so  $G$ is isomorphic to the graph depicted  in Figure \ref{fig3}.
Assume next that   $|G_2|\geq 6$. Then $e(G_2)\leq 18(|G_2|-2)/7$ by the induction hypothesis.  Note that $e(G_1)\leq 3|G_1|-6$ when $|G_1|\leq 5$ and  $e(G_1)\leq 18(|G_1|-2)/7$
 when $|G_1|\ge 6$ by the induction hypothesis.  Therefore, when $|G_1|\le 5$,
\begin{align*}
e(G)&=e(G_1)+e(G_2)\leq 3|G_1|-6+\frac{18(n+1-|G_1|-2)}7\\
&=\frac{18(n-2)}7-\frac{(24-3|G_1|)}7<\frac{18(n-2)}7;
\end{align*}
when $|G_1|\ge 6$,
  \[
  e(G)=e(G_1)+e(G_2)\leq \frac{18(|G_1|+|G_2|-4)}7<\frac{18(n-2)}7.
  \]
So we may assume that  $G$  is $2$-connected. By Observation~\ref{fact}(d), each face in $G$ is bounded by a cycle.
Clearly, $G$ is not a plane triangulation and so $\sum_{i\geq 4}f_i\ge 0$.  We next show that $\sum_{i\geq 5}f_i\ne0$.  Suppose $\sum_{i\geq 5}f_i=0$.  Then $f_3+f_4=f$ and $f_4>0$.  Note that any two distinct $4$-faces of $G$ cannot have exactly one or two edges in common, because $\delta(G)\ge3$ and $G$ is $\Theta_6$-free.  It follows that every edge of a $4$-face of $G$ belongs to $E_{3,4}$, and so $G$ contains a $\Theta_6$ subgraph, a contradiction.
Thus $\sum_{i\geq 5}f_i\ne0$. We may further assume that the  outer face of $G$ is neither a $3$-face nor  a $4$-face. Then
\begin{align*}
2e(G)&=3f_3+4f_4+5f_5+\sum_{i\geq 6}if_i\\
&\geq3f_3+4f_4+5f_5+6(f-f_3-f_4-f_5)\\
&=6f-3f_3-2f_4-f_5,\tag{5}
\end{align*}
which implies that $6f\leq 2e(G)+3f_3+2f_4+f_5$.\medskip

\begin{figure}[htbp]
\centering
\includegraphics [scale=0.55] {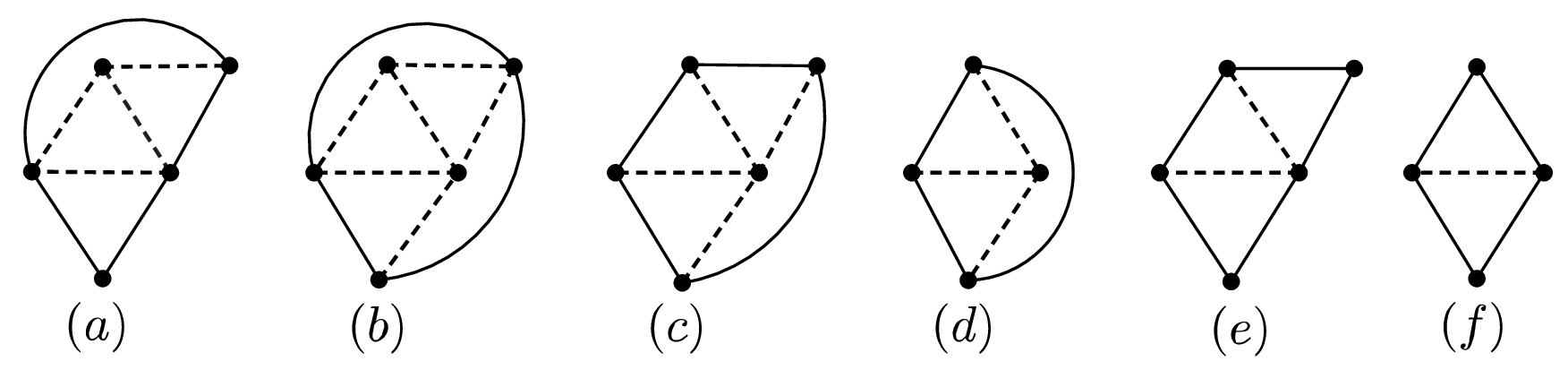}
\caption{All possible configurations of $H_{_F}$, where all dashed edges are in $E_{3,3}$, and   no solid edges are in $E_{3,3}$.}\label{CN}
\end{figure}

We next find an upper bound for each of $f_3, f_4$ and $f_5$. To get an upper bound for $f_3$, we first show that $5e_{3,3}\leq 6f_3$.
 Let $F$ be  a  $3$-face   in $G$ with $|E(F)\cap E_{3,3}|\ge1$. Clearly,   $|E(F)\cap E_{3,3}|\le3$.
Since $G$ is $\Theta_6$-free and the outer face of $G$ is not a $3$-face, there exists  an induced plane subgraph $H_{_F}$ of $G$ with $|H_{_F}|\le5$ such that  $F$ is a face (not outer face) of $H_{_F}$ and all faces of  $H_{_F}$, except the outer face of $H_{_F}$,  are $3$-faces, and no edges on the boundary of the outer face of $H_{_F}$ are  in $E_{3,3}$.  The possible configurations of $H_{_F}$ are shown in Figure~\ref{CN}.
When $H_{_F}$ is isomorphic to  the graph  depicted  in Figure~\ref{CN}(b), $H_{_F}$ contains six edges in $E_{3,3}$ and five $3$-faces of $G$. From all possible  configurations of $H_{_F}$, we see that   $e_{3,3}\leq 6f_3/5$. Hence,
\begin{equation}\tag{6}
3f_3=e_3+e_{3,3}\leq e_3+6f_3/5,  \text{ and so } f_3\leq 5e_3/9.
\end {equation}

\begin{figure}[htbp]
\centering
\includegraphics*[scale=0.25]{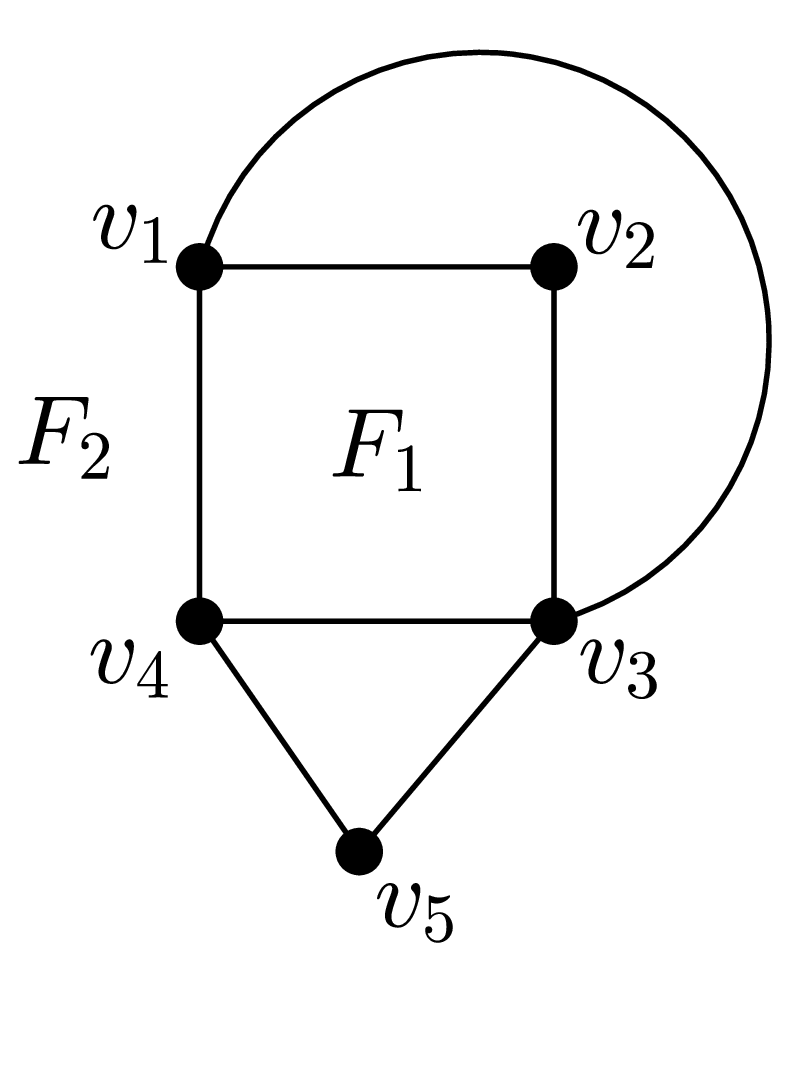}
\caption{Graph $H$.}\label{H}
\end{figure}
To get an upper bound for $f_4$, we first prove  that $e_{4,4}=0$.  Suppose $e_{4,4}>0$.  Let $v_1v_4\in E_{4,4}$,  and let $F_1$ and $F_2$  be two $4$-faces of $G$  having $v_1v_4$  in common.   Then $G$ contains the plane graph $H$, depicted in  Figure~\ref{H},   as a subgraph.   But then $F_2$ is the outer face of $G$, contrary to the assumption that the  outer face of $G$ is neither a $3$-face nor  a $4$-face. Hence, $e_{4,4}=0$.
By Observation~\ref{fact}(b,c),
\begin{equation}\tag{7}
4f_4=e_4 \text{ and }  e(G)\ge e_3+e_4-e_{3,4}.
\end {equation}
We next show that $e_{3,4}\leq e(G)-e_3$. \medskip

This is trivially true when $e_{3,4}=0$.
Assume that  $e_{3,4}\neq0$. Let $F$ and $F'$ be a $4$-face  and a $3$-face in $G$, respectively, such that $F$ and $F'$ share an edge in common. We may assume that $F$ has  vertices $v_1, v_2, v_3, v_4$ in order and $F'$ has  vertices $v_1, v_4, v_5$ in order. Note that $F$ and $F'$ are not outer face in $G$.
Observe that if $v_iv_{i+1}$ belongs to $E_{3,4}$ for any $i\in\{1,2,3\}$, then $v_iv_{i+1}$ belongs to the $4$-face $F$ and the $3$-face with vertices $v_i, v_{i+1}, v_5$ in order,   else $G$ would not be $\Theta_6$-free. Since $n\ge7$, there exists some $k\in \{1,2,3\}$ such that $v_kv_{k+1}\notin E_{3,4}$.  Then $v_kv_{k+1}\in E_{4, j}$ for some $j\ge5$ because $e_{4,4}=0$.  We next show that $F$ has at most two edges in
$E_{3,4}$. Suppose $|E(F)\cap E_{3,4}|=3$. We may assume that $k=2$. Then $v_1v_2,v_3v_4\in E_{3,4}$.  Thus $v_1v_{2}$ belongs to the $4$-face $F$ and the $3$-face with vertices $v_1, v_{2}, v_5$ in order; and $v_3v_{4}$ belongs to the $4$-face $F$ and the $3$-face with vertices $v_3, v_{4}, v_5$ in order. Since $G$ is   $\Theta_6$-free, we see that
 $v_5v_2\in E_{3, j}$ for some $j\ge6$, $v_5v_3\in E_{3, \ell}$ for some $\ell\ge6$, and $v_2v_3\in E_{4, p}$ for some $p\ge6$. But then  $G+v_2v_4$ is $\Theta_6$-free, contrary to the fact that  $G$ is an  extremal  plane graph  for $\Theta_6$ and $n\geq6$. Thus  $F$ has at most two edges in $E_{3,4}$, and so  $F$  has at least two edges not in $E_3$. This holds for each $4$-face in $G$. Hence,    $e_{3,4}\le e(G)-e_3$.\medskip

By (7),
\begin{equation}\tag{8}
 4f_4=e_4\leq e(G)-e_3+e_{3,4}\le 2(e(G)-e_3).
\end{equation}
\medskip
Finally, since $G$ is $\Theta_6$-free, no   $5$-face can share an edge with  a $3$-face in $G$.  Thus $e_{3,5}=0$. By Observation~\ref{fact}(a),
 $e_{5,5}\le e_5\le e(G)-e_3$. By Observation~\ref{fact}(b),
\begin{align*}
5f_5=e_5+e_{5,5}
\leq 2(e(G)-e_3).\tag{9}
\end{align*}
\medskip
Combining $e_3\leq e(G)$ with the upper bounds on $f_3, f_4, f_5$ given in (6), (8), (9), we have
\begin{align*}
6f&\leq 2e(G)+3f_3+2f_4+f_5\\
&\le 2e(G)+5e_3/3+(e(G)-e_3)+2(e(G)-e_3)/5\\
&=17e(G)/5+4e_3/15\\
& \le 11e(G)/3.
\end{align*}
It follows that  $f\leq 11e(G)/18$. By Euler's formula, $n-2\geq e(G)-f\geq 7e(G)/18$. Hence $e(G)\leq 18(n-2)/7$, as desired.\qed
\end{pf}

\begin{cor}\label{strict}  Let  $K_5^-$ be   the graph obtained from $K_5$ by deleting one edge.  Then
\begin{enumerate}[(a)]
\item $ex_{_\mathcal{P}}(n, \Theta_6\cup \{K_5^-\})\leq 12(n-2)/5$  for all $n\ge6$.\vspace{-2mm}
\item $ex_{_\mathcal{P}}(n,\Theta_6)< 18(n-2)/7$ for all $n\geq 10$.
\end{enumerate}
\end{cor}
\begin{pf} To prove (a),  let $G$ be an extremal plane graph for $n\ge6$ and $\Theta_6\cup \{K_5^-\}$. Similar to the proof of Theorem \ref{Theta6}, we see that $e_{3,3}\leq f_3$,  because $G$ is $K_5^-$-free and so no $H_{_F}$ is isomorphic to the graph depicted in Figure~\ref{CN}(b). By Observation~\ref{fact}(b),  $f_3\leq e_3/2$. This, together with the upper bounds for  $f_4, f_5$ given in (8), (9), implies that
\begin{align*}
6f&\leq 2e(G)+3f_3+2f_4+f_5\\
&\le 2e(G)+3e_3/2+e(G)-e_3+2(e(G)-e_3)/5\\
&=17e(G)/5+e_3/10 \\
&\le 7e(G)/2.
\end{align*}
It follows that  $f\leq 7e(G)/12$. By Euler's formula, $n-2\geq e(G)-f\geq 5e(G)/12$. Hence $e(G)\leq 12(n-2)/5$. \\

\begin{figure}[htbp]
\centering
\includegraphics[scale=0.3] {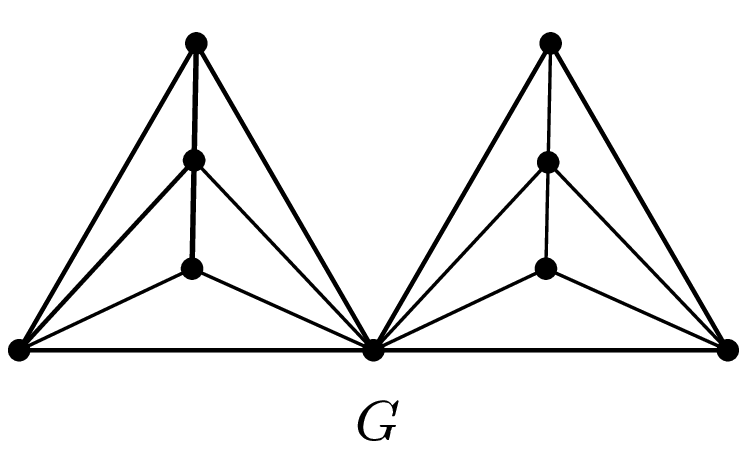}
\hskip 0.5cm
\includegraphics[scale=0.3] {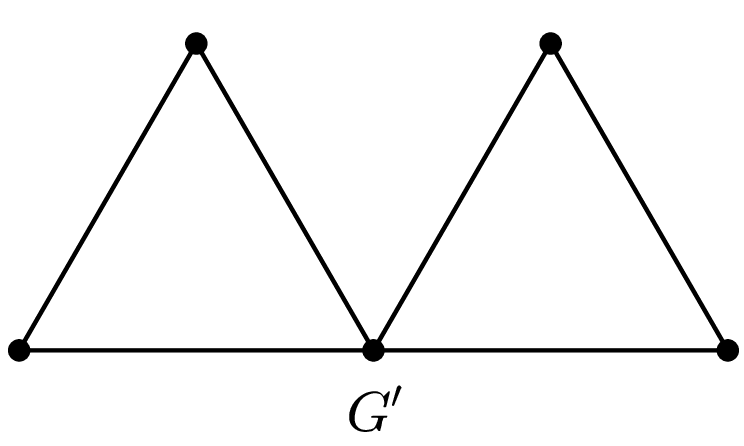}
\caption{An example of constructing the graph $G'$ from a graph $G$.}\label{EX}
\end{figure}

To prove (b), let $G$ be an extremal  plane graph for $\Theta_6$  and  $n\geq 10$. Suppose $e(G)=18(n-2)/7$. By Corollary~\ref{strict}(a), $G$ is not $K_5^-$-free.
From the proof of Theorem~\ref{Theta6}, we see that  equality in $e(G)\le 18(n-2)/7$ is achieved for $n$ if and only if  all the equalities hold in (5), (6), (8), (9) and in $e_3\le e(G)$.  This implies that $e(G)= 18(n-2)/7$ for $n$  if and only if $G$ is a $2$-connected $\Theta_6$-free plane graph on $n$ vertices satisfying: $G$ consists entirely of $K_5^-$'s and $6$-faces,  no two  $K^-_5$'s share an edge,  and no two $6$-faces  have an edge in common. Let $G'$ be the graph obtained from $G$ by deleting the two  vertices not on the outer face in each $H_{_F}=K_5^-$, an example is shown in Figure~\ref{EX}. Then $G'$ consists of $3$-faces and $6$-faces  such that each edge of $G'$ belongs to one $3$-face and one $6$-face. Clearly, $G'$ is $2$-connected because  $G$ is $2$-connected.  Let $f_i'$ be the number of $i$-faces in $G'$. Let $f'=\sum_{i\geq1} f_i'$. Then
$3f_3'=e(G')=6f_6'$ and $f'=f_3'+f_6'.$ Thus $|G'|-2=e(G')-f'=e(G')/2$ and so $e(G')=2|G'|-4$, which implies that $\delta(G')\le 3$. Since $G'$ is $2$-connected, we have  $\delta(G')\ge 2$.  Note that each vertex of $G'$ must have even degree because the adjacent faces are alternatively of size $3$ and size $6$. Thus $\delta(G')= 2$.  Let $v\in V(G')$ be a vertex of degree two in $G'$. Let $u_1vu_2$ and $u_1vu_2u_3u_4u_5$ be the vertices in order  on the boundary of the two adjacent faces containing $v$, respectively.  Then   $G'[\{u_1,v, u_2, u_3, u_4, u_5\}]$ contains a graph in $\Theta_6$ as a subgraph. Thus $G$ is not $\Theta_6$-free, a contradiction.   \qed\medskip
\end{pf}

It is worth noting that   every  $C_6$-free graph  is certainly $\Theta_6$-free. Hence,  $ex_{_\mathcal{P}}(n,C_6)\leq ex_{_\mathcal{P}}(n,\Theta_6)$.   Corollary~\ref{C6} follows immediately from Theorem~\ref{Theta6}.

\begin{cor}\label{C6}
  $ex_{_\mathcal{P}}(n,C_6)\leq{18(n-2)}/7$ for all $n\geq6$, with  equality  when $n=9$.
\end{cor}
\bigskip


\noindent{\bf Acknowledgments.}   \medskip

\noindent  Zi-Xia Song would like to thank Yongtang Shi and  the Chern Institute of Mathematics at Nankai University for hospitality and support during her visit  in May 2017. \medskip

 Yongxin Lan and Yongtang Shi are partially supported by the National Natural Science Foundation of China and Natural Science Foundation of Tianjin (No. 17JCQNJC00300).

\end{document}